\documentclass[conference]{worldcomp}

\usepackage[hmargin=.75in,vmargin=1in]{geometry}
\usepackage[american]{babel}
\usepackage[T1]{fontenc}
\usepackage{times}
\usepackage{caption}
\usepackage{epsfig}

\usepackage{textcomp}
\usepackage{epsfig,graphicx}
\usepackage{xcolor}
\usepackage{amsfonts,amsmath,amssymb}
\usepackage{fixltx2e} 
\usepackage{booktabs}

\title{\bf Greedy Approach for Low-Rank matrix recovery}           

\author{
{\bfseries A. Petukhov$^1$, I. Kozlov$^2$}\\
$^1$Contact Author, Department of Mathematics, University of Georgia, Athens, GA 30602, USA,\\ petukhov@math.uga.edu\\
$^2$Algosoft Tech USA, Bishop, GA, USA, inna@algosoft-tech.com\\
for {\bf IPCV'13}
}

\begin{document}

\maketitle                        

\begin{abstract}
We describe the Simple Greedy Matrix Completion Algorithm providing an efficient method 
for restoration of low-rank matrices from incomplete corrupted entries. 
\par
We provide numerical evidences that, even in the simplest implementation, the greedy approach may 
increase the recovery capability of existing algorithms significantly.
\end{abstract}

\vspace{1em}


\section{Introduction}
We consider a greedy strategy  based algorithm for the recovery of the low-rank
matrix from incomplete corrupted samples.
\par
The problem of low-rank matrix completion is not new.  However,  it got a new impulse  ({\cite{CT2}, \cite{CR}) in connection with the development of the compressed sensing theory and algorithms and ideas to use the ~$\ell^1$ minimization as a surrogate for the~sparsest solution  (\cite{CT}, \cite{D}, \cite{RV}) .  
\par
This paper can be considered as a feasibility study for the methods inspired by ideas from both low-rank
matrix completion and our compressed sensing oriented $\ell^1$-greedy algorithm (\cite{KP1}, \cite{PK2}, \cite{PK3}).
\par
The problem is set as follows. It is required to restore (complete) the matrix $A\in \mathbb R^{m\times n}$   of rank $r$, $r<\min\{m,n\}$,  given by
its $k$ entries, $k<nm$. The set of the given entries is $\Omega\subset \{1,\dots, m\}\times \{1,\dots, n\}$. $|\Omega|$ is the cardinality of $\Omega$ (which
 in our case is equal to $k$). We also introduce notation $d(\Omega)$ for the density of the set $\Omega$, $d(\Omega):=|\Omega|/(nm)$. The complementary set $\bar\Omega$ is a set of erasures, $1-d(\Omega)$ represents the density of erasures.
The theoretical bounds for recoverability of the matrix depends not only on the density of the samples but also on the matrix $A$ and on the 2D-geometry of
$\Omega$. The matrix consisting of only one non-zero entry is the simplest example of a rank 1 matrix which can be restored only if the value at the non-zero entry is known.
Anyway, it turned out (cf. \cite{CR}) that under quite mild conditions random matrices of size $n\times n$ and rank $r$ can be recovered from at most $O(rn^{1.2}\log n)$ entries as a matrix with the minimum of nuclear norm.
\par
The popularity of that 
problem can be explained by an enormous number of applied problems which can be formulated in terms of matrix completion. Among many  settings in different applied areas, we mention problems related to image processing. Image inpainting, including more particular image upsampling, face recognition technique, motion tracking and segmentation in video are most typical of those problems.
While the problem of low-rank matrix completion was studied long ago,  the theory got a big push due to development of Compressed Sensing / Compressive Sampling (CS) technique. After some  simplification,  the CS data decoding goal  can be reduced to solving  underdetermined systems
\begin{equation}
\label{CS}
A\mathbf x=\mathbf y+\mathbf e,
\end{equation}
where $\mathbf x\in\mathbb R^n$ is a sparse vector of data "encoded" with the known to the~decoder matrix  $A\in\mathbb R^{m\times n}$, $m <n$,  $ \mathbf y\in\mathbb R^m$ is a
vector of measurements (of $\mathbf x$) possibly corrupted by the vector $\mathbf e\in\mathbb R^n$. Here and bellow we assume that the sparse solution $\mathbf x$
exists and the vector of errors $\mathbf e$ is also sparse. The sparsity of $\mathbf a\in\mathbb R^N$ means that 
$$
|\mathbf a|_0:= |\{a_i\ne0\}|<N.
$$
The value $|\mathbf a|_0$ is called the Hamming weight of the vector $\mathbf a$.
Since the problem of finding sparse solutions has non-polynomial complexity (\cite{N}), the mainstream CS researches suggested to use to replace the minimization of $|\mathbf x|_0$  (or $|\mathbf x|_0$+$|\mathbf e|_0$) with the minimization of  $\ell^1$-norm. It turned out that that such approach based on convex optimization gives the optimal sparse solution at least when $|\mathbf x|_0$ is
not very large (cf. \cite{D}, \cite{CT}, \cite{RV} for the case $\mathbf e=\vec 0$). Thus, in some special cases the original non-convex problem can be reduced to convex programming. 
In what follows, like for the notion for  the Hamming norm, we use notation $|\cdot|_p$, $0<p<\infty$, for element-wise (quasi-)norms of vectors an matrices. Say, for the matrix $A$, $|A|_p:=\left(\sum_{i,j}|A_{ij}|^p\right)^{1/p}$. In particular, $|\cdot |_2$ is the Frobenius norm. 
The inner product of 2 matrices $A$ and $B$ is defined as $\langle A,B\rangle:=\text{trace}(A^TB)$. Thus, $\langle A,A\rangle=|A|^2_2$.
The notation $\|\cdot\|_p$ is reserved for the operator norms of  matrices.
\par
CS results inspired the authors of \cite{CR} and \cite{CT2} on  replacing  the minimum rank condition leading to non-polynomial complexity with 
the minimization of the nuclear norm $\|A\|_*:=\sum \sigma_i$ of the matrix $A$, where $\sigma_i$ are singular values of $A$. To be more precise, the problem 
$$
\|A\|_*\to \min \text{ subject to } A_{ij}=M_{ij}, \quad (i,j)\in\Omega,
$$
where $M_{i,j}$ are the known entries (measurements) of the matrix $A$, is considered as relaxation of the rank minimization problem above.
\par

Many different settings giving a solution of the original problems have been studied for the last
years. In most cases, the intention of those studies was to find the faster algorithms with the higher
capability of the recovery. Typically, modifications of the problem leading to unconstrained optimization were introduced.
\par
\section{Basic Algorithm}
 For our experiments 
we need the algorithm providing convex minimization recovering low-rank matrices from incomplete corrupted samples. It is used as a basic constructive block in our algorithm. The problem of restoring matrix from corrupted entries is less studied than the simpler matrix completion problem when the available entries are  not corrupted. Anyway, there are a few computationally efficient algorithms solving that problem (e.g.,  \cite{CJSC}, \cite{LCM}, \cite{YYO}).
\par

For our purposes, we selected the algorithm from \cite{LCM} based on the method of Augmented Lagrange Multipliers (ALM)  (e.g., \cite{B}).
Having corrupted samples, instead of finding the matrix with the sparsest set of singular values coinciding with the measurements on as large as possible set, the algorithm finds the minimum of the functional 
$$
L(A,E,Y,\mu):=\|A\|_*+\lambda|E|_1+\langle Y,R\rangle+\frac{\mu}2|R|^2_2,
$$
where $R=M-A-E$ is the residual of approximation of the measurements $M$ of the estimate unknown matrix matrix $A$ and the estimate of the unknown matrix of errors $E$. The entries of the input matrix $M$ on the $\bar\Omega$ are unknown. It is assumed that $R$ vanishes on  $\bar\Omega$ and does not contribute into the  third and the forth term as well as $E$ does not contribute into the second term.
\par
If it is known that the observed entries in $M$ are not corrupted, the second term can be omitted.
However, we assume that we never know whether the entries are corrupted. So, in what follows, we 
minimize the 4-term functional. 
\par
We will need the following notation
$$
\mathcal S_{\epsilon}[x]:=\left\{\begin{array}{ll} x-\epsilon, &x>\epsilon,\\x+\epsilon, &x<-\epsilon,\\0, &\text{otherwise};\end{array}\right.
$$
where $x$ can be either a number or a vector or a matrix. The operator $\mathcal S_{\epsilon}$ is called the shrinkage operator. The norms $\|\cdot\|_2$, $\|\cdot\|_{\infty}$ applied to matrices mean the operator norms.
\par
The minimization algorithm, as it is described in  \cite{LCM} and implemented in Matlab code, is as follows

{\bf Algorithm ALM.}
\par
{\bf Input.} Observation matrix $M\in \mathbb R^{m\times n}$, defined on $\Omega$, and $\lambda>0$.
\par
{\bf Initialization.} $Y^0=\frac{1}{\max\{\|M\|_2,\|M\|_{\infty}/\lambda \}} M$, \par $E^0=0$, $\mu_0>0$, $\rho>1$, $k=0$;
\par
1. {\bf while} not converged {\bf do}
\par
2. $(U,S,V):=\text{svd}(M-E^k+\mu^{-1}Y^k)$;
\par
3. $A^{k+1}:=U\mathcal S_{\mu_k^{-1}}[S]V^T$;
\par
4.  $E^{k+1}:=\mathcal S_{\lambda\mu_k^{-1}}[M-A^{k+1}+\mu^{-1}Y^k]$;
\par
5. $Y^{k+1}:=Y^k+\mu_k(M-A^{k+1}-E^{k+1})$;
\par
6. $\mu_{k+1}:=\rho\mu_k$, $k:=k+1$;
\par
7. {\bf end while}.
\par
{\bf Output.} $A^{k+1}$, $E^{k+1}$.
\par
\section{Our algorithm}
Our modification of the algorithm above is inspired by significant success reached by applying
greedy ideas to solving underdetermined systems (\cite{KP1}, \cite{PK2}, \cite{PK3}). The general greedy strategy in optimization 
algorithms consists in sequential finding a simple suboptimal solutions giving some information  about the optimal solution. A greedy algorithm picks up some most obvious features or elements
of those solutions and gives  them a privilege to be pivot for the next iteration of the suboptimal algorithm. Each iteration brings new pivot elements. 
\par
In the matrix recovery algorithm, the erasures from the set $\bar\Omega$ forms such group from the beginning. Whereas, the elements of $\Omega$ are just suspicious to be erroneous.  If we have  sufficient evidence that some element in  $\Omega$  contains a random error independent of the content of other entries from $\Omega$, that the decision to move this element to $\bar\Omega$ is 
quite justifiable. While the independence condition is not always accurate even in our experiments
with artificially generated data, we use this strategy for estimation of the capability of the greedy ideas for matrix recovery. 
\par
Our greedy algorithm consists in iterating  with an updated (dilated) sets $\Omega_k$. We will call it the Simple Greedy Matrix Completion  Algorithm (SGMCA). Generally speaking , any matrix recovery algorithm, including SGMCA  itself, which is able to fight the mixture of erasures and errors can be used as a basic block of SGMCA. 
\par
Formally, all our experiments can be described in the following way.

{\bf Algorithm SGMCA.}

{\bf Input.} $M$, $\Omega$.
{\bf Initialization.}  $\lambda>1$, $\Omega_0:=\Omega$, $A^0:=M$, $E^0:=0$, $0<q<1$, $k=0$.

1. Set $k:=k+1$;

2.  $A^k=ALM(M,\Omega_k)$;

3.  {\bf if} $k=1$  {\bf then} $T_1=0.3\max_{i,j}\{|A^1_{ij}-M_{ij}|$; {\bf else} $T_k:=0.65T_{k-1}$;

 4. $\Omega_{k+1}=\Omega_k \setminus \{(i,j)\mid |A^k_{ij}-M_{ij}|>T_k\}$;

5. {\bf if} not converged {\bf go to} 2.

\section{Numerical Experiments}
Since our intention was to conduct algorithm feasibility  study, the goal of this section is to give comparison with the output of recently published algorithms and with pure ALM (one iteration of the algorithm above with no update).  
The parameters in the basic algorithm are selected as $\mu_0=0.3/\|M\|_2$, $\rho=1.1+0.5|\Omega|/(mn)$. 
The parameter $\lambda$ is defined by the combination of $d(\Omega)$ and the density of errors in $\Omega$ samples. The general trend can be characterized as follows: the  higher error rate, the less value of $\lambda$ has to be used. In what follows, we do not use fine tuning of $\lambda$. 
The same $\lambda$ is used for big groups of experiments. At the same time, tuning  $\lambda$ may bring significant increase of the algorithm efficiency. In this paper, we use values of $\lambda$ in the range $0.02 \div 100$
\par
For our experiments, we used Matlab implementation of ALM algorithm available at 
are available at { \tt http://perception.csl.illinois.edu} {\tt/matrix-rank/home.html}
We used the code for Matrix Completion via the inexact ALM Method with our adaptation to the input with errors.
\par 
In all our experiments, $A$ are square $m\times n$, $m=n$ matrices of rank $r$ obtained 
as $A=UV^T$, where $U,V\in\mathbb R^{n\times r}$. The matrices $U$ and $V$ consist of 
independent gaussian random values with zero expectation and the variance 1.
The coordinates of erasures were selected randomly. The models of errors below were different for different experiments.

In the first experiment (Fig.1), we demonstrate advantage of the iterative SGMCA over ALM (one iteration of the same algorithm). For the matrix with fixed sizes $m\times n$, $m=n$, and the rank of matrices 
$r=15$. The solid curves on Fig. 1  correspond to  SGMCA (up to 10 iterations) for  $n=128,\,512,\, 1024$ (from the bottom to the top). The corresponding graphs for ALM algorithm are plotted with dashed curves. 

\noindent
\begin{picture}(00,210)
\label{Fig1}
\put(0,10){          \epsfig{file=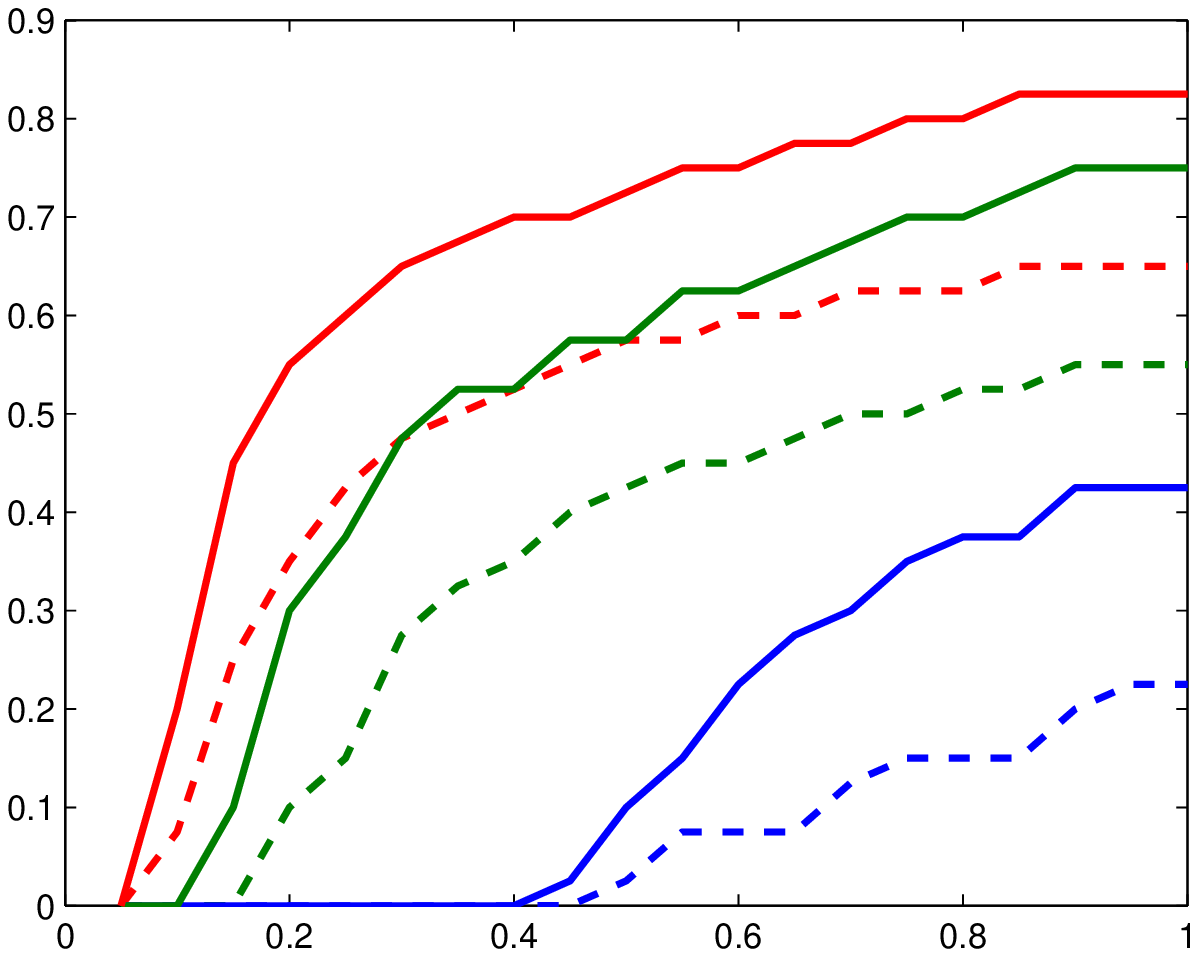,
           height=200pt, width=260pt}
		  } 
\put(50,10){\footnotesize {\bf Fig.1. Greedy iterations vs. convex ALM optimization}}
\end{picture}
\par
The horizontal coordinate indicates the fraction of the matrix available for restoration (i.e., $d(\Omega)$), while the vertical
coordinate is the fraction of randomly  corrupted  entries  in $\Omega$. The magnitude of the corruption is randomly set from the standard normal distribution. The curves define "phase transition" bounds. In our experiments, we run 10 trials. At the points of curves
all 10 attempts were accomplished with success, i.e., for the obtained estimate $\hat A$, $|A-\hat A|_2/|A|_2<10^{-3}$, whereas for the points above the curves, at least one attempt failed. It is assumed that the regions under the curves are regions of "success".
\par
The second experiment (Fig.2 and Fig.3) is devoted to comparison with the results from \cite{CJSC}. 
Unfortunately, we do not have full information about the error model. So we use the same  additive model of errors as above. While all other parameters are taken from \cite{CJSC}.
The matrix of rank 2 is constructed as above.  Its size changes from 100 to 3000.

The experiment  consists of 2 parts. The first plot (Fig.2) contains the curves for the  fixed erasure rate 0.1, i.e, $d(\Omega)=0.9$. However, the probability of errors in those entries varies.  There are 3 graphs on Fig. 2. The solid line corresponds to 10 iterations of SGMCA, the dashed  line corresponds to ALM, and the dotted curve corresponds to the result from \cite{CJSC}. The values defined by curves give the  maximum error probability admitting successful correction by the corresponding algorithm. If we were aware of the error model   from   \cite{CJSC}, the dashed and dotted curves have to coincide up to statistical discrepancy. 
\par
On the second plot (Fig.3) the error rate is fixed and equal to 0.1. The graphs show dependence of maximum possible rate of erasures from the size of matrix. 

\par
\noindent
\begin{picture}(00,210)
\label{Fig2}
\put(0,10){          \epsfig{file=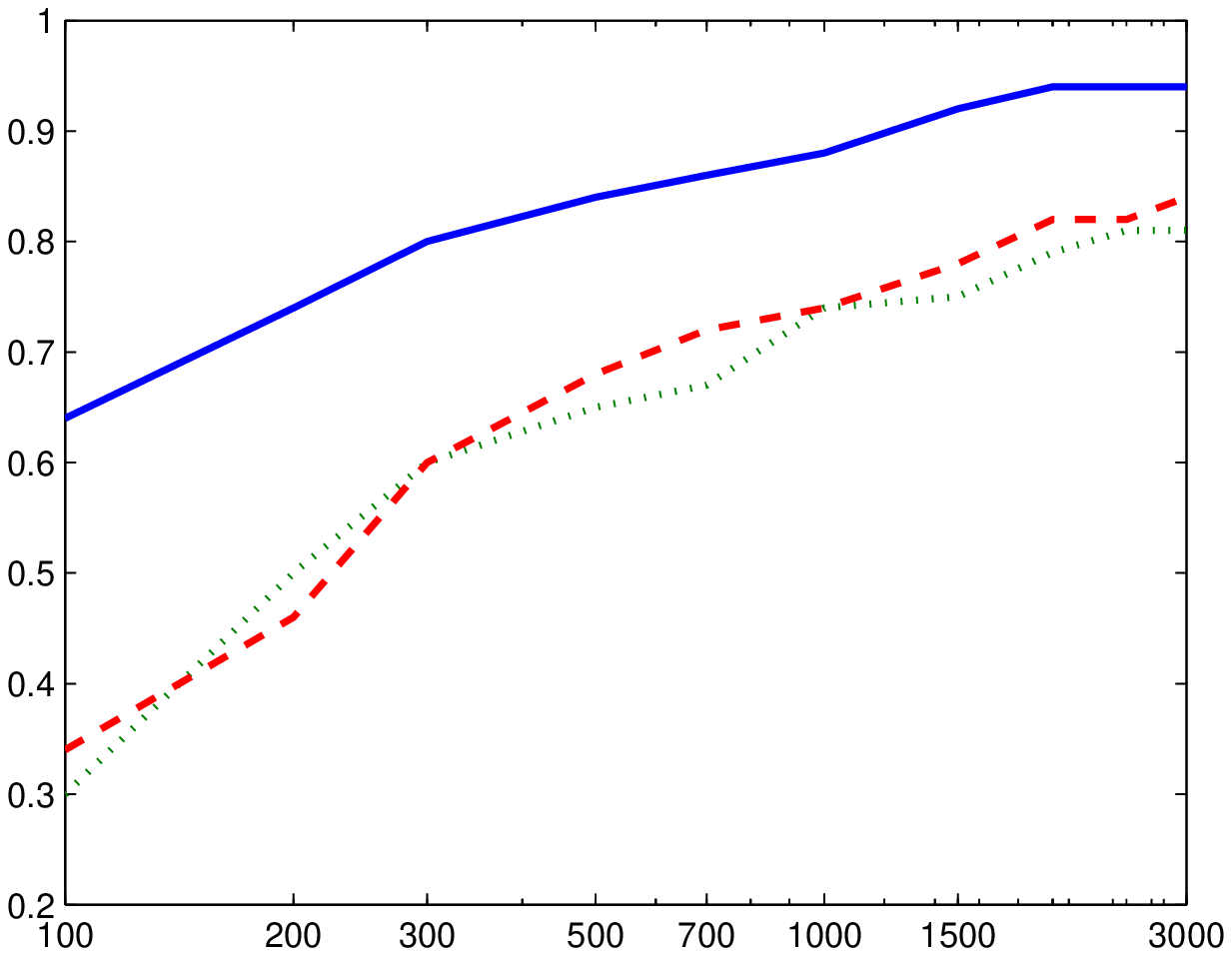,
           height=200pt, width=260pt}
		  } 
\put(50, 5){\footnotesize {\bf Fig.2. Admissible error rate for erasure rate 0.1.}}
\end{picture}
\par
\par
 The efficiency of the algorithms is defined by the distance of  curve values on  Fig. 2 and Fig. 3 to 1. It is easy to see that, when the error rate is fixed, SGMCA curve is twice closer to 1 than ALM. Hence, SGMCA  restores low rank matrices  from only half of entries necessary for the  ALM minimization.  For the fixes erasure rate, the number of uncorrupted entries for SGMCA successful restoration can be only  one third of that necessary for recovery with ALM.

\par
\noindent
\begin{picture}(00,210)
\label{Fig2}
\put(0,10){          \epsfig{file=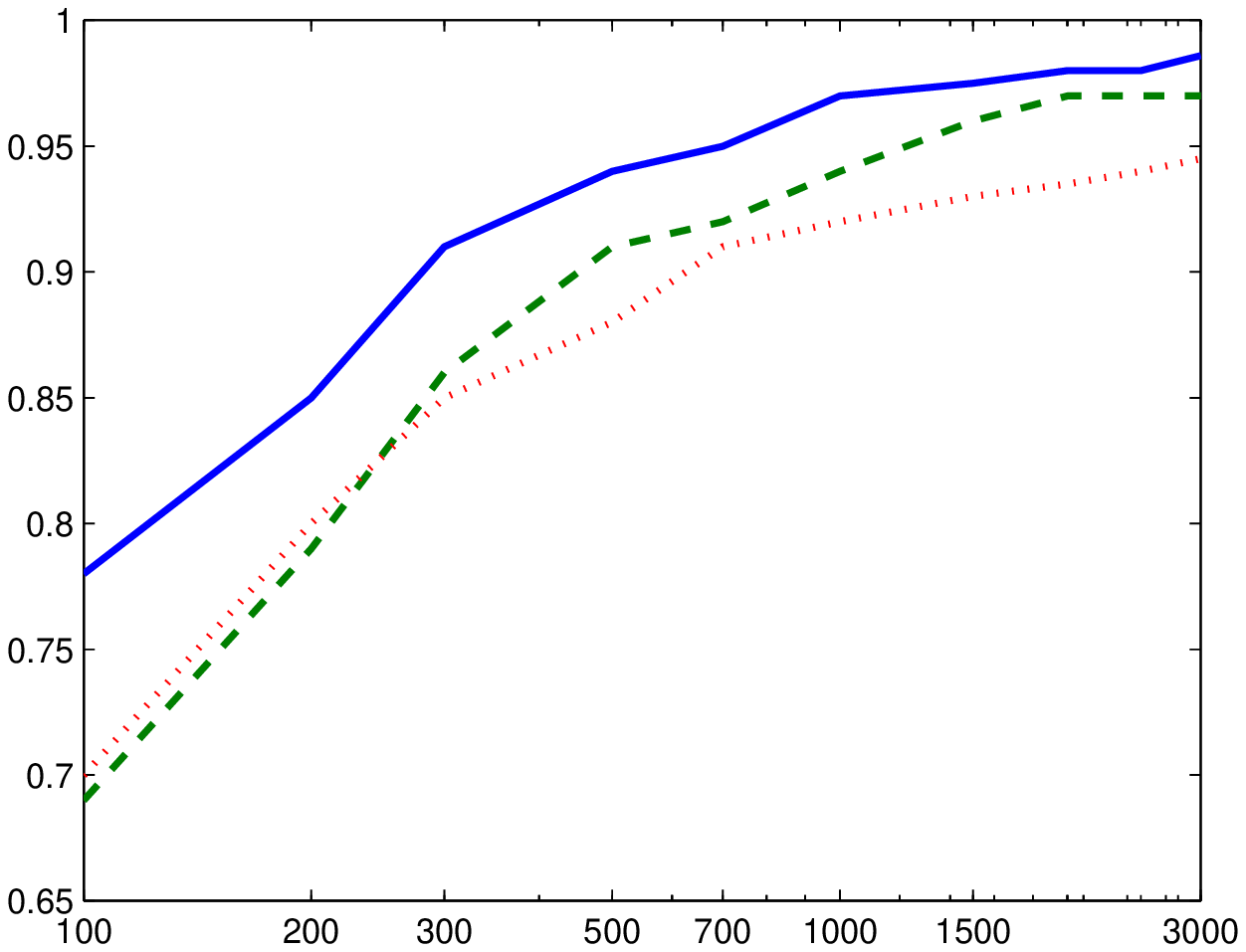,
           height=200pt, width=260pt}
		  } 

\put(50, 5){\footnotesize {\bf Fig.3. Admissible erasure  rate for error rate 0.1.}}

\end{picture}

\par
The more precise value of the SGMCA graph at $m=3000$  on Fig 3 is 0.986, i.e., having error probability 0.1, the matrix can be restored from 1.4\% of random entries.
\par
In our last experiment (Fig. 4--6), we compare the output of SGMCA with the results of RTRMC algorithm from \cite{YYO} 
 providing very impressive recovery. However, for such successful
recovery it requires  a priori knowledge of the rank of the matrix $A$. The ALM-based algorithm used as a basic algorithm in SGMCA does not require any knowledge about the rank while the rank knowledge is useful for it. To provide  equal opportunities we applied the  internal fixation of the~rank inside ALM procedure. 
\par
The results for ranks $r= 5,\,15,\,  25$ are given on Fig.  4--6     correspondingly. The size of matrices is $512\times512$. The horizontal coordinate is $d(\Omega)$, whereas the the vertical coordinate is the probability of errors in the coefficients available for reconstruction. In most of cases, SGMCA outperforms RTRMC by 15--25\% in the maximum admissible probability of errors. The reason why SGMCA loses on interval $[0, 0.175]$ on Fig.4 is the parameter $\lambda=0.02$ which was fixed for all 3 experiments. Setting $\lambda=0.2$ on that interval, we would get the overwhelming advantage of SGMCA. We emphasize that when the rank is known in advance,  the optimal parameter $\lambda$ can be computed for each $d(\Omega)$. This would not contradict the equal opportunity of the~two algorithms. Optimal selection of $\lambda$ is a reserve not used in our experiments.
\par
The model of data in  \cite{YYO}  is identical to the model described above.  Whereas, the error model is different.  The values of the corrupted values are randomly uniformly distributed between minimum and maximum of uncorrupted values. We also use that error model in our experiment.  
\par
\noindent
\begin{picture}(00,210)
\label{Fig2}
\put(0,10){          \epsfig{file=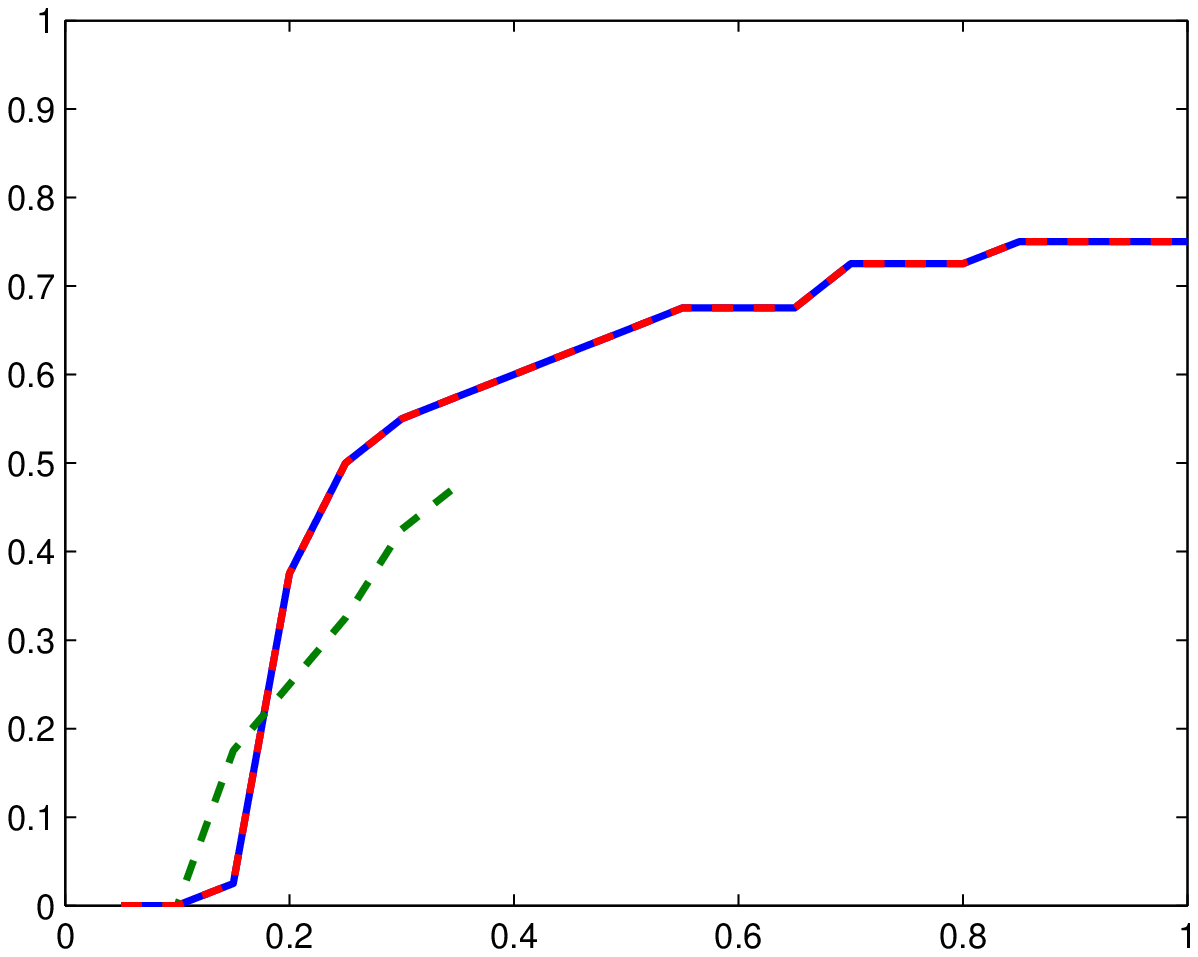,
           height=200pt, width=260pt}
		  } 

\put(50, 5){\footnotesize {\bf Fig.4. SGMCA vs. RTRMC. $r=5$.}}

\end{picture}

\par
\noindent
\begin{picture}(00,210)
\label{Fig2}
\put(0,10){          \epsfig{file=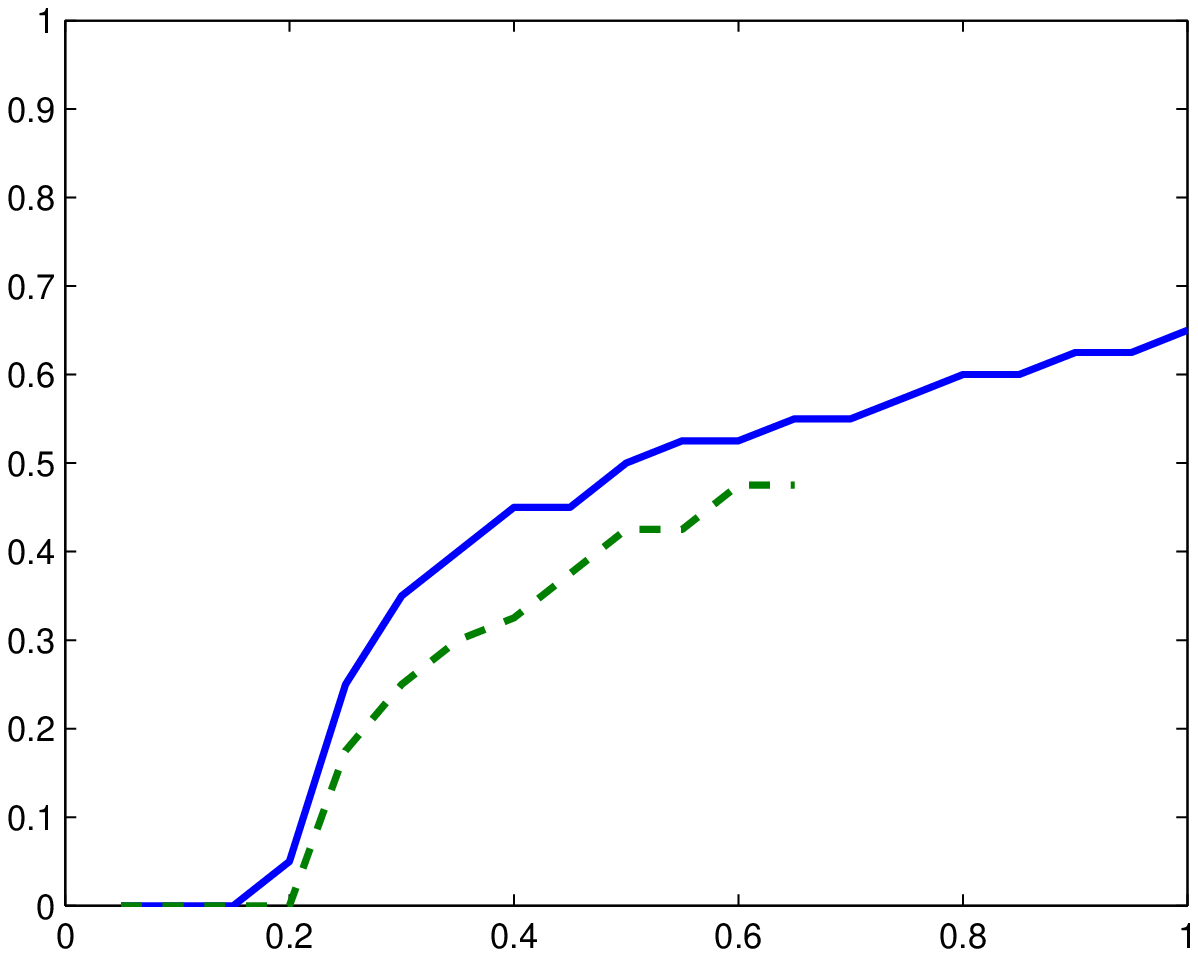,
           height=200pt, width=260pt}
		  } 

\put(50, 5){\footnotesize {\bf Fig.5. SGMCA vs. RTRMC. $r=15$.}}

\end{picture}

\par
\noindent
\begin{picture}(00,210)
\label{Fig2}
\put(0,10){          \epsfig{file=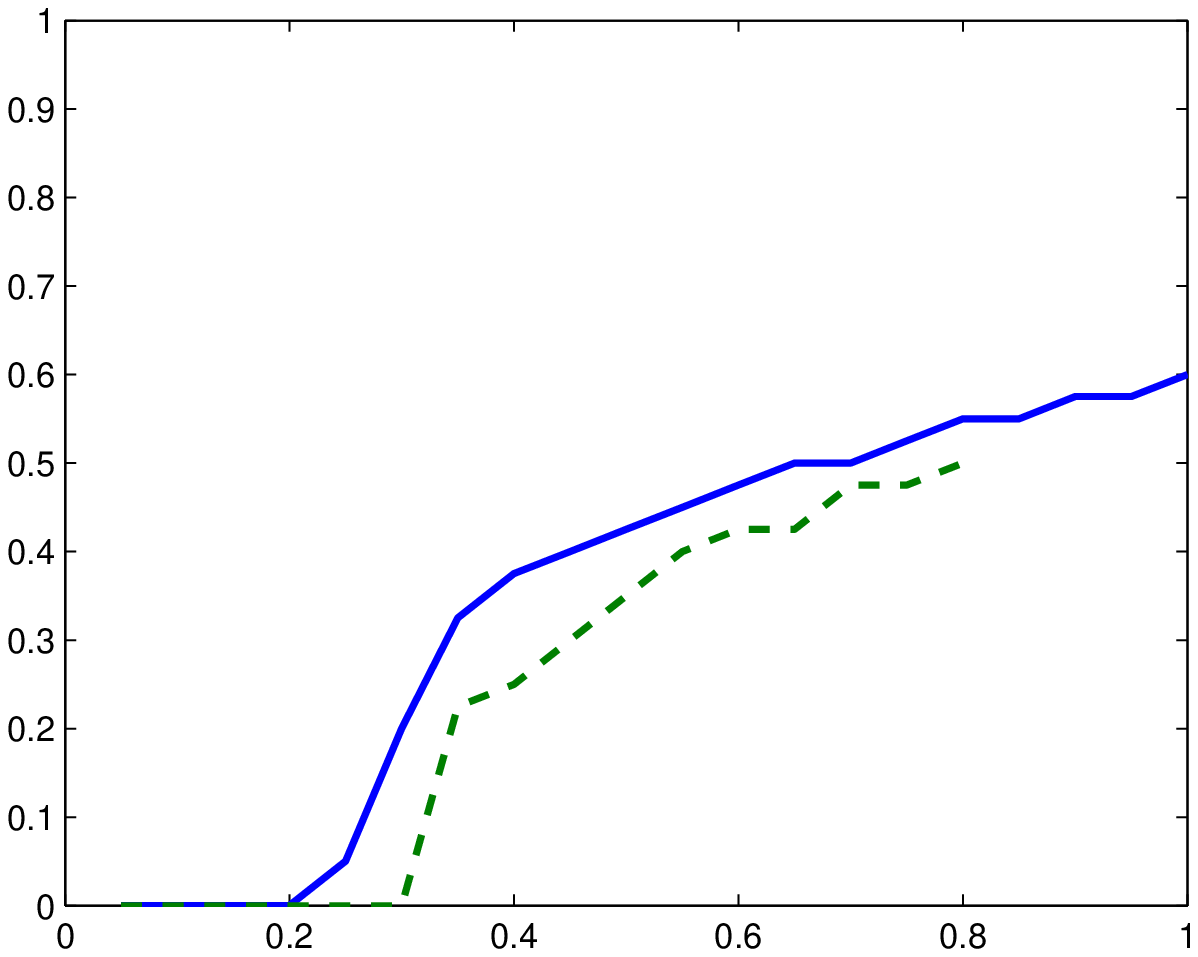,
           height=200pt, width=260pt}
		  } 

\put(50, 5){\footnotesize {\bf Fig.6. SGMCA vs. RTRMC. $r=25$.}}

\end{picture}

\par
The dashed line corresponding to RTRMC is shorter than our solid line since we used the data directly from \cite{YYO}.
\par
\section{Future Studies}
This study shows that even the simplest implementation of the greedy idea in the
 form of SGMCA outperforms the recent algorithms significantly in the recovering ability for very incomplete measurements with high level of corruption. The results above show the feasibility of the idea, its perspectives and a high level of expectation. 
\par
Because of its iterative nature, the algorithm has to repeat basic step a few times. We restricted a number of iterations by 10 however in most of cases 5 iterations provided necessary precision. 
Among possible directions for improvement, acceleration ways are needed to be considered.
One of possible ways is to do not wait for the completion of each iteration, updating $\Omega$ within 
internal iterations of the basic algorithm. In the simplest, but maybe in less efficient form, it can be done even with no  intrusion into the basic algorithm.  For instance, when the greedy step looks for coordinates 
of large errors, we do not need high precision output of the basic algorithm. So  an update of the precision on each iteration from low to high may accelerate the algorithm for the data close to 
the limits of the potential recovering ability (phase transition points). However, it should be noted that this 
modification may slightly slow down the recovery of the data located far from the phase transition points. Other way for acceleration skipped by us  in this paper is to use the estimate of $A$ on the previous iteration as a basic algorithm start point  for the next iteration. 
\par
Now we discuss the ways for increasing the capability SGMCA in the matrix recovery.
\par
First of all, in our experiments, we practically did not use fine tuning of the weight $\lambda$. We used fixed $\lambda$ for big range of the parameters of input data. Whereas, just by replacing the value 0.02 with the value $0.2$, the results presented on Fig.4--6 can be significantly improved for the high level of erasures (on the left side of the graphs). Indeed, when the rank of $A$ and the model of errors is known, the optimal values of $\lambda$ admitting the maximum density of error can be found in advance for any number of erasures and used in the recovery procedure. 
\par
When the rank is unknown or nature of possible corruption is unknown in advance, adaptive finding $\lambda$ becomes a challenging problem. This problem has many common features with the problem of finding sparse solutions of  underdetermined linear systems with corrupted input. In the mentioned problem, the weight $\lambda$ is defined by
the interaction between the sparsity of the possible solution and the error vector. In the problem of 
of matrix completion, the low rank plays a role of the solution sparsity   in CS. So the methods (or at least principles) developed in CS can be applied to the matrix completion problem. In \cite{PK3}, we showed that the sparsity of the solution can be reliably estimated on the dynamic of change 
of the value $|\mathbf x|_{0.5}/|\mathbf x|_{2}$. The same characteristic of the matrix $A$ intermediate estimates can be used for the matrix completion procedure. We also have to say that, generally speaking, for finding $\lambda$ we do not need both the rank and the sparsity of errors estimates. Indeed, if we have the sparsity of errors, we can evaluate the potential maximum  rank $r$ admitting recovery with the given algorithm provided that $\lambda$ is optimal.  Then that optimal $\lambda$ will also provide recovery of matrices with the rank  less than $r$. Thus, the adaptive $\lambda$ is one of quite realistic sources for increasing algorithm capability.
\par
\section{Conclusion}
The paper presents feasibility study for the Simple Greedy Matrix Completion algorithm. We showed that it outperforms recently developed algorithms of matrix completion significantly. We also discussed the ways for further increase SGMCA efficiency.

\end{document}